\documentclass[12pt]{article}
\usepackage{amsmath,amssymb,amsthm}
\usepackage{amsfonts, amsmath}
\newtheorem{theorem}{Theorem}[section]
\newtheorem{lemma}{Lemma}[section]
\newtheorem{corollary}{Corollary}[section]

\newtheorem{proposition}{Proposition}[section]
\theoremstyle{definition}
\newtheorem{definition}{Definition}[section]
\theoremstyle{remark}
\newtheorem{remark}{Remarks}[section]
\newtheorem{example}{Example}[section]
\numberwithin{equation}{section}

\title{On summability, integrability  and impulsive differential equations in Banach spaces}
\date{}
\author {S. Heikkil\"a\\
Department of Mathematical Sciences, University of Oulu\\
BOX 3000, FIN-90014, Oulu, Finland\\
E-mail: sheikki@cc.oulu.fi}

\begin{document}
\maketitle 

\begin{abstract} {\bf Purpose:} To study summability of families indexed by  well ordered sets of $\mathbb R\cup\{\infty\}$ in normed spaces. To derive integrability criteria for step mappings and for right regulated mappings from an interval of $\mathbb R\cup\{\infty\}$ to a Banach space. To study solvability of impulsive differential equations.
\newline
{\bf Main methods:} A generalized iteration method presented, e.g., in \cite{HL94}. Summability of families  in normed spaces indexed with well ordered subsets of $\mathbb R\cup\{\infty\}$.
\newline
{\bf Results:} Necessary and sufficient conditions for global and local HK, HL, Bochner and Riemann integrability of step mappings and right  regulated mappings defined on an interval of $\mathbb R\cup\{\infty\}$, and having values in a Banach space. Applications to impulsive differential equations are also presented. Families  indexed with well ordered subsets of $\mathbb R\cup\{\infty\}$ are used to represent  impulsive parts of considered equations and to approximate their solutions.  

\noindent {{\bf MSC:} 26A06, 26A18, 26A39, 26B12, 26E20, 34A36, 34A37, 34G20, 40A05, 40D05, 40F05, 47H07, 47H10}

\noindent{\bf Keywords:} family; well ordered index set; summable; primitive; integral;
 step mapping; right regulated; differential equation; impulsive; generalized iteration method. 

\noindent{\bf Acknowledgements:} The author thanks a Referee for very valuable comments.

\end{abstract}

\vskip3pt


\newpage

\baselineskip 18pt

\section{Introduction}\label{1} In Chapter VIII of his book "Foundations of Modern Analysis" Jean Dieudonn$\acute{e}$  criticised the adoption of the  Riemann integral in Calculus courses  as follows: "Only the stubborn conservatism of academic tradition could freeze it into a regular part of curriculum, long after it had outlived its historical importance." 
The integral presented in the book is justified as follows:  
"to dispense with any considerations drawn from measure theory; this is what we have done by defining only the integral of regulated mappings (sometimes called the "Cauchy integral"). When one needs a more powerful tool, there is no point in stopping halfway, and the general theory of ("Lebesgue") integral is the only sensible answer."

On the other hand, a few years before the publication  in 1960 of the book \cite{Die60} cited above,  Ralph Henstock and Jaroslav Kurzweil generalized the definition of the
Riemann integral so that the resulting integral, called here the Henstock--Kurzweil (shortly HK) integral, encloses the Lebesgue integral of real valued functions. The study of HK integrals of  Banach valued mappings started around 1990 by the work of R. A. Gordon. The strong version of HK integral is called here the Henstock--Lebesgue (shortly HL) integral.
No measure theory is needed in the definitions of HK and HL integrals. Moreover, if a mapping $g$ from a compact real interval $I$ to a Banach space $E$ is Bochner integrable, i.e., if the norm function $t\mapsto \|g(t)\|$ is Lebesgue integrable, then
$g$ is also HL integrable. Converse is not true because the norm function of a HL integrable mapping is not necessarily  HL integrable or HK integrable. Moreover HL integrability encloses improper integrals on finite intervals; HK integrability also on unbounded intervals.
 
In \cite{Die60} the integral calculus is presented for regulated mappings, i.e., the mappings from a real interval $I$ to a Banach space $E$, having left limits in $I\setminus\{\inf I\}$ and right limits in $I\setminus\{\sup I\}$. The definition of the integral of a mapping  $g:I\to E$ is based on the existence of a primitive, i.e., a continuous mapping $f:I\to E$ that is differentiable in the complement of a countable subset $Z$ of $I$, and $f'(t)=g(t)$,  for each $t\in I\setminus Z$.  Because 
any two primitives of $g$ differ by a constant, the difference $f(b)-f(a)$ for any two points of $I$, is independent of the particular primitive $f$. This difference is written $\int_a^bg(t)\,dt$, and is called the integral of $g$ from $a$ to $b$. As shown in \cite{Die60}, a primitive exists for every  regulated mapping. 

In this paper we study integrability of  right regulated mappings, i.e., those mappings from an interval $I$ of $\mathbb R\cup\{\infty\}$ to a Banach space $E$ which have right limits at every point of $I\setminus\{\sup I\}$. The main difference between regulated mappings and right regulated mappings is that the latter ones may have discontinuities of the second kind, while regulated mappings can have only discontinuities of the first kind. 
Another difference is that regulated mappings are  HL integrable and Riemann integrable on bounded intervals, whereas all right regulated mappings are not even HK integrable.  The main purpose of this paper is to develop criteria for HK, HL, Bochner and Riemann integrability of right regulated mappings on an interval $I$ of $\mathbb R\cup\{\infty\}$.
Necessary and sufficient conditions for local integrability 
are also presented. 
The main tools are:
\begin{itemize}
\item   A generalized iteration method presented, e.g., in \cite{HL94}. Using this method we shall prove that a right regulated mapping has at most countable number of discontinuities, and that it can be approximated uniformly on compact intervals by step mappings with well ordered steps. A fixed point theorem based on this method is applied in the study of impulsive differential equations.
\item   Summability of families in normed spaces. Summability of families with nonempty index sets is studied, e.g, in \cite{Bour1,Bour2,Die60}.  The given definitions rule out conditional summability, so that the obtained summability results are not applicable in the study of HK and HL integrability.
Therefore we assume that the index set is well ordered. For the sake of applications we assume that
the index set is contained in $\mathbb R\cup\{\infty\}$.
\item CD primitives. By a CD {\it primitive} of a mapping $g$ from an interval $I$ of $\mathbb R\cup\{\infty\}$ to $E$ we mean  a continuous mapping $f:I\setminus\sup I\to\mathbb R$ that is differentiable in the complement of a countable subset $Z$ of $I$, and $f'(t)=g(t)$  for each $t\in I\setminus Z$.
\end{itemize}

This paper is organized as follows. In Section \ref{2} we define and study  summability and absolute summability of a family $(x_\alpha)_{\alpha\in\Lambda}$ in a normed space  when the index set $\Lambda$ is a well ordered subset of $\mathbb R\cup\{\infty\}$. 
With the help of such families   we derive necessary and sufficient conditions for global and local HK, HL, Bochner and Riemann integrability of step mappings and right  regulated mappings defined on an interval of $\mathbb R\cup\{\infty\}$, and having values in a Banach space. The results obtained for step mappings in Section \ref{S3} both generalize and improve some results derived in \cite{DPM02,Sye05,S04} (see Remark \ref{R301}). The integrability criteria derived in Section \ref{S4} for right regulated mappings are new. We shall prove, for instance, the following results for a right regulated mapping $g:I\to E$, $-\infty<\min I<\sup I\le\infty$. (We say that a property holds locally for a function defined on $I$, if the function has that property on every compact subinterval of $I$.) 
\begin{enumerate}
\item  $g$ is locally HL integrable if and only if it has a CD primitive.
\item  $g$ is HL integrable when $I$ is bounded if and only if $g$ has a CD primitive that has the left limit at $\sup I$.
\item  $g$ is HK integrable if it has a CD primitive that has the left limit at $\sup I$.
\item  $g$ is locally Bochner integrable if and only if it has a locally absolutely continuous CD primitive.
\item  $g$ is Bochner integrable if and only if it has a locally absolutely continuous CD primitive that has the left limit at $\sup I$.
\item  $g$ is locally Riemann integrable if and only if it is locally bounded, in which case $g$ has a locally Lipschitz continuous CD primitive.
\item $g$ is Riemann integrable if and only if it is bounded and $I$ is bounded.
\item The improper Riemann integral of $g$ from $\min I$ to $\sup I$ exists if  $g$ is locally bounded, and its CD primitive has the left limit at $\sup I$.  
\item For each compact subinterval $[a,b]$ of $I$, either $g$ is Riemann integrable on $[a,b]$, or  there exists the greatest number $c_1$ in $(a,b]$ such that $g$ is locally Riemann integrable on $[a,c_1)$.
\item For each compact subinterval $[a,b]$ of $I$, either $g$ is Bochner integrable on $[a,b]$, or   there exists the greatest number $c_2$ in $(a,b]$ such that $g$ is locally Bochner integrable on $[a,c_2)$.
\item For each compact subinterval $[a,b]$ of $I$, either $g$ is HL integrable on $[a,b]$, or there  exists the greatest number $c_3$ in $(a,b]$ such that $g$ is 
locally HL integrable on $[a,c_3)$. 
\end{enumerate}
%

Concrete examples of mappings $f,\,g:\mathbb R_+\to E$ are presented for  above results when
$E$ is the space $c_0$ of those sequences of real numbers which converge to $0$. In every example the mapping $g$ has the discontinuity of second kind at every rational point of $\mathbb R_+$.  
The above results are valid with minor modifications also when $g$ is left regulated, i.e., 
when $g$ has left limits at every point of $I\setminus\{\inf I\}$. 

The first one of the above results will be applied in Section \ref{5}
to impulsive  differential equations
in Banach spaces. Families  indexed with well ordered subsets of $\mathbb R\cup\{\infty\}$ are used to represent  impulsive parts of considered equations and to approximate their solutions.     

\section{Preliminaries}\label{2}
\setcounter{equation}{0}
In this section we shall first present basic properties of well ordered subsets of $\mathbb R\cup\{\infty\}$. These sets are used as index sets of families in normed spaces. After defining summability and presenting examples of such families  we introduce basic facts on HK, HL, Bochner and Riemann integrability of mappings from a  real interval to a Banach space.
  
A nonempty subset $\Lambda$ of $\mathbb R\cup\{\infty\}$,  ordered by the natural ordering $<$ of $\mathbb R$, plus  $t<\infty$ for every $t\in\mathbb R$, is well ordered if    
every nonempty subset of $\Lambda$ has the smallest element. In particular, to every number $\beta$ of 
$\Lambda$, different from its possible maximum, there corresponds the smallest element in $\Lambda$ that is greater than $\beta$. It is called the successor of $\beta$ and is denoted by $S(\beta)$.  
There are no numbers of $\Lambda$ in the open interval $(\beta,S(\beta))$.  
The following properties are needed:
\begin{itemize}
\item Every well ordered subset of $\mathbb R\cup\{\infty\}$ is countable.

\item Principle of Transfinite Induction: If $\Lambda$ is well ordered and $\mathcal P$ is a property such that if  $\mathcal P(\gamma)$ is true whenever $\mathcal P(\beta)$ is true for all $\beta < \gamma$ in $\Lambda$, then
$\mathcal P(\gamma)$ is true of all $\gamma\in \Lambda$.
\end{itemize}

\begin{definition}\label{D1} Let $\Lambda$ be a well ordered 
subset of $\mathbb R\cup\{\infty\}$. Denote $a=\min \Lambda$, and  $b=\sup\Lambda$. 
When $\gamma\in\Lambda\cup\{b\}$, denote $\Lambda^{<\gamma}=\{\alpha\in\Lambda|\alpha < \gamma\}$. 
The family $(x_\alpha)_{\alpha\in\Lambda}$ with elements $x_\alpha$ in a normed space $E$  is summable if it has the following properties: 
\begin{itemize}
\item[(s)] To every $\gamma\in\Lambda\cup\{b\}$ there corresponds a unique element $\sigma(\gamma)$ of $E$, called the sum of the family $(x_\alpha)_{\alpha\in\Lambda^{<\gamma}}$, satisfying the following conditions:
\item[(i)] $\sigma(a)=0$, and if $\gamma=S(\beta)$ for some $\beta\in\Lambda$, then 
$\sigma(\gamma)=\sigma(\beta)+x_\beta$.
\item[(ii)] If $\gamma$ is not a successor, then for each $\epsilon > 0$  there is such  $\beta_\epsilon\in\Lambda^{<\gamma}$ that $\|\sigma(\beta) -\sigma(\gamma)\|
<\epsilon$ whenever $\beta\in\Lambda$ and $\beta_\epsilon\le \beta<\gamma$.
\item[] The sum $\sigma$ of a summable family $(x_\alpha)_{\alpha\in\Lambda}$ is $\sigma(b)$ if $b\not\in\Lambda$, and $\sigma(b)+x_b$ if $b\in\Lambda$.
\item[] If $\sigma(\gamma)$ is defined for every $\gamma\in \Lambda$, we say that $(x_\alpha)_{\alpha\in\Lambda}$ is locally summable.
\item[] A family $(x_\alpha)_{\alpha\in\Lambda}$ is (locally) absolutely summable if  $(\|x_\alpha\|)_{\alpha\in\Lambda}$ is (locally) summable.
\end{itemize}
\end{definition}

\begin{remark}\label{R}
The above definition of summability is analogous to that given in \cite{Wiki} when the index set $\Lambda$ is an ordinal.
Because $\Lambda$ is countable, the given definition of absolute summability of  a family $(x_\alpha)_{\alpha\in\Lambda}$
 agrees on that of \cite[Section V.3]{Die60}, i.e., for a bijection $\varphi$ from $\mathbb N$ to $\Lambda$ the series $\overset{\infty}{\underset{n=1}{\sum}}x_{\varphi(n)}$ is absolutely convergent. 
As for results dealing with  ordinary, unconditional and absolute convergence of $\overset{\infty}{\underset{n=0}{\sum}} z_n$ when $E$ is a Banach space, see, e.g., \cite[Appendix B]{Sye05}.
\end{remark} 

Next we shall determine the first partial sums and the sum of a summable  family $(x_\alpha)_{\alpha\in\Lambda}$ in some elementary cases ($\sup_\Lambda$ means the least upper bound in $\Lambda$).

\begin{enumerate}
\item If $\Lambda$ is finite and nonempty, then $\Lambda=\{S^n(a)|n=0,\dots,m\}$, $m\in\mathbb N_0$ ($S^0(a)=a$). 
\item If $\Lambda=\{S^n(a)|n\in\mathbb N_0\}$, then $\sigma=\overset{\infty}{\underset{n=0}{\sum}} x_{S^n(a)}=\overset{\infty}{\underset{n=0}{\sum}} z_n$, where $z_n=x_{S^n(a)}$, $n\in\mathbb N_0$.
\item After  $S^n(a)$, $n\in\mathbb N_0$, the next possible numbers of $\Lambda$ are $a_0=\sup_\Lambda\{S^n(a)|n\in\mathbb N_0\}$, $S^m(a_0)$, $m=0,1,\dots$, $a_1=\sup_\Lambda\{S^m(a_0)|n\in\mathbb N_0\}$, $\dots$, $S^m(a_i)$, $m=0,1,\dots$, $a_{i+1}=\sup_\Lambda\{S^m(a_i)|n\in\mathbb N_0\}$, $i=0,1,\dots$, $b_0=\sup_\Lambda \{a_i|i\in\mathbb N_0\}$, and so on. Corresponding 
partial sums of the family  $(x_\alpha)_{\alpha\in\Lambda}$ are:
$$
\begin{aligned}
&\sigma(S^m(a_0))=\sigma(a_0)+ \sum_{n=0}^{m-1}x_{S^n(a_0)}= \sum_{n=0}^\infty x_{S^n(a)}+ \sum_{n=0}^{m-1}x_{S^n(a_0)}, \\ 
&\sigma(a_1)= \sigma(a_0)+ \lim_{m\to\infty}\sum_{n=0}^{m-1}x_{S^n(a_0)}= \sum_{n=0}^\infty x_{S^n(a)}+ \sum_{n=0}^{\infty}x_{S^n(a_0)}, \\    
&\sigma(a_{i+1})= \sigma(a_i)+ \lim_{m\to\infty}\sum_{n=0}^{m-1}x_{S^n(a_i)}= \sum_{n=0}^\infty x_{S^n(a)}+ \sum_{j=0}^i\left(\sum_{n=0}^{\infty}x_{S^n(a_j)}\right), \\ 
&\sigma(b_0)= \sum_{i=0}^\infty\sigma(a_i)=\sum_{n=0}^\infty x_{S^n(a)}+ \sum_{j=0}^\infty\left(\sum_{n=0}^{\infty}x_{S^n(a_j)}\right),
\end{aligned}
$$
and so on. In particular, if $b=\sup\Lambda=b_0$, we have the associative rule: $\sigma(b)=\overset{\infty}{\underset{i=0}{\sum}}\sigma(a_i)$, where the
sum of $(x_\alpha)_{\alpha\in\Lambda^{<b}}$ is presented as a sum of an infinite number of its partial sums.
However, this presentation is not independent on the order of both partial sums and their elements, as in the case of absolutely or unconditionally summable families.   
\end{enumerate}

\begin{example}\label{Ex201} A simple example of a well ordered subset of an interval $[a,b)$ of $\mathbb R$ is an increasing sequence formed by numbers
\begin{equation}\label{E201}
b-2^{-n}(b-a), \quad n\in\mathbb N_0.
\end{equation}
The smallest number of this sequence is $a$ and its supremum is $b$. When $a=0$ and $b=1$ we obtain the sequence
$$
\Lambda_0=\{\alpha(n_0)=1-2^{-n_0}|n_0\in\mathbb N_0\}.
$$
The points of $\Lambda_0$ divide the interval $[0,1)$ into disjoint subintervals $[1-2^{-n_0},1-2^{-n_0-1})$, $n_0\in\mathbb N_0$. Choosing
$a=1-2^{-n_0}$, $b=1-2^{-n_0-1}$ and $n=n_1$ in (\ref{E201}) we obtain in each of these subintervals decreasing sequences,
which together form an inversely well ordered set  
$$
\Lambda_1=\{\alpha(n_0,n_1)=1-2^{-n_0-1}-2^{-n_0-n_1-1}|n_0,n_1\in\mathbb N_0\}.
$$
Choosing a vector $e\ne 0$ of $E$ and denoting
$$
x_{\alpha(n_0,n_1)}=\frac{(-1)^{n_1}2^{-n_0}e}{n_1}, \quad n_0,n_1\in\mathbb N_0,
$$
we obtain a summable  family
$$
(x_{\alpha(n_0,n_1)})_{\alpha(n_0,n_1)\in\Lambda_1}= \sum_{n_0=0}^\infty2^{-n_0}\left(\sum_{n_1=0}^\infty\frac{(-1)^{n_1}e}{n_1}\right).
$$
The above process can be continued in the obvious way. For each $m\in\mathbb N_0$ one obtains a well ordered
set 
\begin{equation}\label{E202}
\Lambda_m=\bigg\{\alpha(n_0,\dots,n_m)=1-\sum_{k=0}^{m-1}2^{-\sum_{j=0}^kn_j-j-1}-2^{-\sum_{j=0}^mn_j-m}\bigg|n_0,\dots,n_m\in\mathbb N_0\bigg\}.
\end{equation}
Denoting
$$
x_{\alpha(n_0,\dots,n_m)}=\frac{(-1)^{n_m}2^{-\sum_{k=0}^{m-1}n_k}e}{\sqrt[m]{}n_m}, \quad n_0,\dots,n_m\in\mathbb N_0,
$$
then the  family
$$
(x_{\alpha(n_0,\dots,n_m)})_{\alpha(n_0,\dots,n_m)\in\Lambda_m}= \sum_{n_0=0}^\infty\left(\cdots\left(\sum_{n_{m-1}=0}^\infty\left(\sum_{n_m=0}^\infty\frac{(-1)^{n_m}2^{-\sum_{k=0}^{m-1}n_k}e}{\sqrt[m]{}n_m}\right)\right)\right)
$$
is summable but neither absolutely nor unconditionally summable.

In above considerations  $\min\Lambda_m=0$ $\sup\Lambda_m=1$ and  for every $m\in\mathbb N_0$. Another way is to
restrict $\Lambda_0$ to $[0,\frac 12)$, $\Lambda_1$ to $[\frac 12,\frac 34)$, and in general, restrict $\Lambda_m$ to $[1-\frac 1{2^{m}},1-\frac 1{2^{m+1}})$, $m\in\mathbb N_0$. Thus $\alpha(n_0,\dots,n_m)$ is replaced by 
$\beta(n_0,\dots,n_m)=2^{-m-1}(1+\alpha(n_0,\dots,n_m))$, i.e.,
$$
\beta(n_0,\dots,n_m)=1-2^{-m-1}-\sum_{k=0}^{m-1}2^{-\sum_{j=0}^kn_j-j-m-2}+2^{-\sum_{j=0}^mn_j-2m-1}.
$$
These numbers form a well ordered set
$$
\Lambda^1_0=\{\beta(n_0,\dots,n_m)|m,n_0,\dots,n_m\in\mathbb N_0\}
$$
satisfying $\min\Lambda^1_0=0$, and $\sup\Lambda^1_0=1$. 

Replacing in the above considerations $\Lambda_0$ by $\Lambda^1_0$ we obtain more general  well ordered sets of rational numbers: $\Lambda^1_m$, $m\in\mathbb N_0$, $\Lambda^2_0$,\dots, $\Lambda_0^n$, $n\in\mathbb N_0$.
For $n>1$ a  family $(x_\alpha)_{\alpha\in\Lambda_0^n}$ is no more representable as a multiple series.    
\end{example}

The following result is needed in the integrability studies.

\begin{lemma}\label{L0.00} Let $(x_\alpha)_{\alpha\in\Lambda}$ be a family in  $E$ having a well ordered index set $\Lambda$ in $\mathbb R\cup\{\infty\}$. 
\newline
(a) Either $(x_\alpha)_{\alpha\in\Lambda}$ is bounded, or there is the greatest element $c_1$ in $\Lambda\setminus\{\min\Lambda\}$
 such that the family $(x_\alpha)_{\alpha\in\Lambda^{<\gamma}}$ is bounded for each $\gamma\in\Lambda^{<c_1}$. 
\newline
(b) Either $(x_\alpha)_{\alpha\in\Lambda}$ is absolutely summable, or there is the greatest element $c_2$ in $\Lambda\setminus\{\min\Lambda\}$
 such that the family $(x_\alpha)_{\alpha\in\Lambda^{<\gamma}}$ is absolutely summable for each $\gamma\in\Lambda^{<c_2}$.
\newline
(c) Either $(x_\alpha)_{\alpha\in\Lambda}$ is  summable, or there is the greatest element $c_3$ in $\Lambda\setminus\{\min\Lambda\}$
 such that the family $(x_\alpha)_{\alpha\in\Lambda^{<\gamma}}$ is summable for each $\gamma\in\Lambda^{<c_3}$.
\newline
(d) $c_1$, $c_2$ and $c_3$ are not successors.
\end{lemma}    

{\bf Proof}. (a) If $(x_\alpha)_{\alpha\in\Lambda}$ is not bounded, there is at least one number $c$ in $\Lambda$ such that  $(x_\alpha)_{\alpha\in\Lambda^{<c}}$ is not bounded. Because $\Lambda$ is well ordered, there is the smallest of such numbers $c$. Denoting it by $c_1$, then  the family $(x_\alpha)_{\alpha\in\Lambda^{<\gamma}}$ is bounded for each $\gamma\in\Lambda^{<c_1}$, but not for each $\gamma\in\Lambda^{<c}$, if $c_1<c\in\Lambda$. This proves (a).
\newline
(b) Assume that the family $(x_\alpha)_{\alpha\in\Lambda}$ is not absolutely summable. Then there is at least one number $c$ in $\Lambda$ such that  $(x_\alpha)_{\alpha\in\Lambda^{<c}}$ is not absolutely summable. Because $\Lambda$ is well ordered, there is the smallest of such numbers $c$. Denoting it by $c_2$, then  the family $(x_\alpha)_{\alpha\in\Lambda^{<\gamma}}$ is absolutely summable for each $\gamma\in\Lambda^{<c_2}$, but not for each $\gamma\in\Lambda^{<c}$, if $c_2<c\in\Lambda$. This proves (b). 
\newline
(c) The proof of (c) is similar to that of (b) when absolute summability is replaced by summability.
\newline
(d) To prove that $c_1$ is not a successor, assume on the contrary that $c_1=S(c)$ for some $c\in \Lambda$. Thus $(x_\alpha)_{\alpha\in\Lambda^{<c_1}}=(x_\alpha)_{\alpha\in\Lambda^{<c}}\cup\{x_c\}$. Because $(x_\alpha)_{\alpha\in\Lambda^{<c_1}}$ is unbounded, then $(x_\alpha)_{\alpha\in\Lambda^{<c}}$ also unbounded.
But $c < S(c)=c_1$, whence $c_1$ is not the smallest of those numbers $c$ of $\Lambda$ for which $(x_\alpha)_{\alpha\in\Lambda^{<c}}$ is unbounded, contradicting with the choice of $c_1$.

The proofs that $c_2$ and $c_3$ are not successors  are similar.
\qed    

A mapping $g$ from a compact real interval
$[a,b]$ to a Banach space $E$ is called Henstock-Lebesgue
(shortly HL) integrable if there is a mapping $f:[a,b]\to E$, called a
{\em primitive} of $g$, with the following property: To each
$\epsilon > 0$ there corresponds such a mapping $\delta:[a,b]\to
(0,\infty)$ that whenever $[a,b]=\overset{m}{\underset{i=1}{\cup}}[t_{i-1},t_{i}]$  and
$\xi_i\in [t_{i-1},t_i]\subset
(\xi_i-\delta(\xi_i),\xi_i+\delta(\xi_i))$ for all $i=1,\dots, m$,
then
\begin{equation}\label{E000.2}
\sum_{i=1}^m \bigg\|f(t_i)-f(t_{i-1})-g(\xi_i)(t_i-t_{i-1})\bigg\|<\epsilon.
\end{equation}
$g$ is called Henstock-Kurzweil (shortly HK) integrable if the above property holds with (\ref{E000.2}) replaced by
$$
\left\|\sum_i
(f(t_i)-f(t_{i-1})-g(\xi_i)(t_i-t_{i-1}))\right\|<\epsilon.
$$
Primitives of HK and HL integrable mappings are continuous (see \cite[Theorem 7.4.1]{Sye05}). 
If $g$ is   HL  (resp. HK) integrable on $[a,b]$, it is  HL (resp. HK) integrable on every
closed subinterval $[c,d]$ of $[a,b]$. Because 
any two primitives of $g$ differ by a constant, the difference $f(c)-f(c)$ for any two points of $[a,b]$, is independent of the particular primitive $f$. This difference is  called the {\it Henstock--Kurzweil integral} of $g$ from $c$ to $d$, and is denoted by  $\int_c^dg(s)\,ds$. Thus, 
\begin{equation}\label{E200}
\int_c^dg(s)\,ds:=f(d)-f(c),\ \hbox{ where $f$ is a primitive of
$g$.}
\end{equation}
Riemann integrability  is obtained when in the definition of HK integrability the gauge functions $\delta$ are replaced by  positive constants $\delta$. In this case the integral, defined by (\ref{E200}), is called the {\it Riemann integral}.

As for the proofs of the following results, see, e.g., \cite[Proposition 24.44 and Theorem 24.45]{Sch97}. 

\begin{lemma}\label{L000} A mapping $g:[a,b]\to E$ is Riemann integrable if $g$ is bounded, and is continuous in the complement of a subset $Z$ of $[a,b]$ that has Lebesgue measure $0$. Conversely, every Riemann integrable mapping is bounded.
\end{lemma}

A mapping $g:I\to E$, $-\infty< \min I< \sup I\le\infty$, is said to be locally integrable in HK, HL, Bochner or Riemann sense if $g$ is HK, HL, Bochner or Riemann integrable on every compact subinterval of $I$.  

The next lemma follows, e.g., from \cite[Lemma 1.12]{CH11}.

\begin{lemma}\label{L000.0}
If a mapping $g:I\to E$ has a CD primitive $f$, then $g$ is locally HL
integrable, and (\ref{E200}) holds for every compact subinterval $[c,d]$ of $I$.
\end{lemma}

As for the definition of the HK integral on unbounded real intervals, and the proof of the next result, see \cite{BoSa04}.

\begin{lemma}\label{L000.1} If $-\infty< a< b\le\infty$ and $g:[a,b]\to E$, then the following results are equivalent.
\newline
(a) $g$ is HK integrable on $[a,c]$ for each $c\in[a,b)$, and  $\underset{c\to b-}{\lim}\int_a^c g(s)\,ds$ exist.
\newline 
(b)  $g$ is HK integrable on $[a,b]$, and $\int_{a}^bg(s)\,ds=\underset{c\to b-}{\lim}\int_a^c g(s)\,ds$. 
\end{lemma}

\begin{remark}\label{R703.00}
By definition every HL integrable mapping is HK integrable. Converse holds if $E$ is finite dimensional.
(see \cite[Proposition 3.6.6]{Sye05}). If $b<\infty$, the result of Lemma \ref{L000.1} holds when HK integrability is replaced by HL integrability.

A strongly measurable mapping $g:[a,b]\to E$ is Bochner integrable if and only if the function $t\mapsto \|g(t)\|$ is Lebesgue integrable. Every Bochner integrable mapping is HL integrable. 
In particular, HL integrability encompasses improper integrals of
Bochner integrable mappings.

Regulated mappings are HK, HL, Bochner and Riemann integrable.
\end{remark}

In the proof of  the following lemma we apply a generalized iteration method. 

\begin{lemma}\label{L400} Let $g:[a,b]\to E$ be right regulated. Then to every positive number $\epsilon$ there corresponds such a  well ordered set $\Lambda_\epsilon$ in $[a,b]$ that $[a,b)$ is a disjoint union of half-open intervals $[\beta,S(\beta))$, $\beta\in \Lambda_\epsilon^{<b}$, and $\|g(s)-g(t)\|\le \epsilon$ whenever $s,\,t\in (\beta,S(\beta))$ and $\beta\in \Lambda_\epsilon^{<b}$.  
\end{lemma}

{\bf Proof}. Define $G_\epsilon:[a,b]\to [a,b]$ by $G_\epsilon(b)=b$, and
\begin{equation}\label{E402}
G_\epsilon(x) = \sup\{y\in(x,b]|\ \|g(s)-g(t)\|\le \epsilon\ \hbox{ for all } \ s,\,t\in (x,y)\}, \quad x\in[a,b).
\end{equation}
It is easy to verify that $G_\epsilon$ is increasing, i.e., $G_\epsilon(x)\le G_\epsilon(y)$ whenever $a\le x\le y\le b$. Because $g$ is right regulated, then $x<G_\epsilon(x)$ for each $x\in[a,b)$.
By \cite[Theorem 1.1.1]{HL94} there is exactly one well ordered subset $\Lambda_\epsilon$ of $[a,b]$ having the following property:
\begin{equation}\label{E403}
a=\min \Lambda_\epsilon,\ \hbox{ and $a <\gamma\in \Lambda_\epsilon$ if and only if } \ \gamma=\sup \{G_\epsilon[\{\beta\in \Lambda_\epsilon|\beta< \gamma\}]\}.
\end{equation}
Because $\sup G_\epsilon[\Lambda_\epsilon]$ exists, it is by \cite[Theorem 1.1.1]{HL94} both a fixed point of $G_\epsilon$ and $\max \Lambda_\epsilon$. Since $b$ is the only fixed point of $G_\epsilon$, then $b=\max \Lambda_\epsilon$.
Since $\beta < G_\epsilon(\beta)$ for each $\beta\in\Lambda_\epsilon^{<b}$, it follows from \cite[Lemma 1.1.3]{HL94}  that $G_\epsilon(\beta)=S(\beta)$ for all $\beta\in\Lambda_\epsilon^{<b}$. Thus, by  \cite[Corollary 1.1.1]{HL94},  
 $[a,b)$ is the disjoint union of half-open intervals $[\beta,S(\beta))$, $\beta\in \Lambda_\epsilon^{<b}$. The last conclusion follows from (\ref{E402}),
since $G_\epsilon(\beta)=S(\beta)$ for all $\beta\in \Lambda_\epsilon^{<b}$.
\qed

With the help of Lemma \ref{L400} we shall prove some properties for right regulated mappings.

\begin{lemma}\label{L401} Let $g:[a,b]\to E$ be right regulated. Then 
\newline
(a) $g$ has at most a countable number of discontinuities. 
\newline
(b) $g$ is strongly measurable.
\end{lemma}

{\bf Proof}. (a) Let $\Lambda_n$, $n\in\mathbb N$,  denote the well ordered subset $\Lambda_\epsilon$ defined by  (\ref{E403}) when $\epsilon=\frac 1n$.
 It follows from  Lemma \ref{L400} that  $\|g(s)-g(t)\|\le\frac 1n$ whenever  $s,\,t\in(\beta,S(\beta))$ and $\beta\in \Lambda_n^{<b}$. Thus all the discontinuity points of $g$ belong to the countable set $Z=\bigcup_{m=1}^\infty \Lambda_m$.   

(b) By (a) the set $Z$ of discontinuity points of $g$ is a null set, whence $g$ is strongly measurable.
\qed

\section{On HL, HK, Bochner and Riemann integrability of step mappings}\label{S3}
\setcounter{equation}{0}

 Let $E$ be a Banach space. In this section we consider first the integrability of a step mapping $g:[a,b]\to E$, $-\infty <a < b <\infty$, that has well ordered steps, i.e., there is a well ordered subset $\Lambda$ of $[a,b]$ such that
$\min\Lambda=a$ and $\max\Lambda = b$, and a family $(z_\alpha)_{\alpha\in\Lambda}$ of $E$ such that
\begin{equation}\label{E30}
g(t)=z_\alpha,\ t\in[\alpha,S(\alpha)),\quad \alpha\in\Lambda^{<b}.
\end{equation}
Assume also that $[a,b)$ is a countable union of disjoint intervals $[\alpha,S(\alpha))$, $\alpha\in\Lambda$. Thus $g$ is well-defined on $[a,b)$ by (\ref{E30}). 

As an application of Lemma \ref{L000.0} we shall prove the following result.

\begin{proposition}\label{P301} Assume that $g:[a,b]\to E$, $-\infty <a < b <\infty$,  is a  step mapping with representation (\ref{E30}) on $[a,b)$. 
Then the following condition are equivalent:
\newline
(a) $g$ is HL integrable.
\newline
(b) The family $((S(\alpha)-\alpha)z_\alpha)_{\alpha\in\Lambda^{<b}}$ 
is summable.
\newline
If (a) or (b) holds, then $\int_a^bg(t)\,dt$ is the sum of the family $((S(\alpha)-\alpha)z_\alpha)_{\alpha\in\Lambda^{<b}}$. 
\end{proposition}

{\bf Proof}.  Assume first that the family $((S(\alpha)-\alpha)z_\alpha)_{\alpha\in\Lambda^{<b}}$ 
is summable. Denote by 
 $\sigma(\gamma)$ the sum of  $((S(\alpha)-\alpha)z_\alpha)_{\alpha\in\Lambda^{<\gamma}}$, $\gamma\in\Lambda$. 
We shall show that the mapping $f:[a,b)\to E$, defined  by 
\begin{equation}\label{E000.4}
f(t)=\sigma(\gamma)+ (t-\gamma)z_\gamma, \ t\in[\gamma,S(\gamma)), \quad \gamma\in\Lambda^{<b},
\end{equation}
is a CD primitive of $g$. It follows from (\ref{E30}) and (\ref{E000.4}) that 
$$
f'(t)=z_\gamma=g(t), \quad t\in(\gamma,S(\gamma)), \ \gamma\in \Lambda^{<b}.
$$
Thus $f'(t)=g(t)$ for every $t\in [a,b)\setminus\Lambda^{<b}$. In particular, $f$ is continuous in $[a,b)\setminus\Lambda^{<b}$.  
To prove that $f$ is continuous at every point of $[a,b)$, it suffices to prove continuity at every point $\gamma\in\Lambda^{<b}$.
Since $f(t)=(t-a)z_a$, $a\le t < S(a)$, then  $f$ is right continuous at $\gamma = a$. If $\gamma\in\Lambda^{<b}$ is a successor, i.e., $\gamma=S(\beta)$ for some $\beta\in \Lambda^{<b}$, then  
$$
f(t)=\begin{cases}\sigma(\beta)+ (t-\beta)z_\beta, \ t\in[\beta,\gamma)),\\
                 \sigma(\gamma)+ (t-\gamma)z_\gamma, \ t\in[\gamma,S(\gamma)).\end{cases}
$$
Applying the defining condition (s) of summability we obtain
$$
\lim_{t\to\gamma-}f(t)= \sigma(\beta)+ (S(\beta)-\beta)z_\beta =\sigma(\gamma)= \lim_{t\to\gamma+}f(t).
$$
Thus $f$ is continuous at $\gamma=S(\beta)$, $\beta\in\Lambda^{<b}$. 

Assume next that $\gamma$ is not a successor. Given $\epsilon > 0$, there exists by condition (s)(ii) of summability such a $\beta_\epsilon \in\Lambda^{<\gamma}$ that 
$$
\|\sigma(\beta) -\sigma(\gamma)\|
<\epsilon \ \hbox{whenever $\beta\in\Lambda$ and} \ \beta_\epsilon\le \beta<\gamma.
$$
If $t\in(\beta_\epsilon,\gamma)$, there exists $\beta\in \Lambda$, $\beta_\epsilon\le \beta<\gamma$, such that $t\in[\beta,S(\beta))$. Thus
$$
\|f(t)-f(\gamma)\|=\|\sigma(\beta)+(t-\beta)z_\beta-\sigma(\gamma)\|<\epsilon +\|(t-\beta)z_\beta\|.
$$
Since also $\beta_\epsilon\le S(\beta)<\gamma$, and since 
$$\|(t-\beta)z_\beta\|\le \|(S(\beta)-\beta)z_\beta\|=\|\sigma(S(\beta))-\sigma(\beta)\|,
$$
then
$$
\|f(t)-f(\gamma)\|<\epsilon+\|\sigma(S(\beta))-\sigma(\beta)\|\le\epsilon+\|\sigma({S(\beta))}-\sigma(\gamma)\|+\|\sigma(\beta)-\sigma(\gamma)\|<3\epsilon.
$$
This holds for every $t\in[\beta_\epsilon,\gamma)$. Thus $\underset{t\to\gamma-}{\lim}f(t)=f(\gamma)$.
 If $t\in[\gamma,S(\gamma))$, then 
$$
f(t)=\sigma(\gamma)+ (t-\gamma)z_\gamma, \ t\in[\gamma,S(\gamma)).
$$ 
Thus $\underset{t\to\gamma-}{\lim}f(t)=\sigma(\gamma)=f(\gamma)=\underset{t\to\gamma+}{\lim}f(t)$.
This proves that $f$ is continuous at $\gamma$.

The above proof shows that $f$ is continuous in $[a,b)$, and that $f'(t)=g(t)$ in the complement of the well ordered, and hence countable subset $\Lambda^{<b}$ of $[a,b)$. Thus $f$ is a CD primitive of $g$, so that $g$ is locally HL integrable on $[a,b)$ by  Lemma \ref{L000.0}. Using condition (s) it can be shown (cf. the proof of Proposition \ref{P32}) that  $f(t)\to \sigma(b)$ as $t\to b-$. 
Thus $\int_a^tg(s)\,ds=f(t)-f(a)=f(t)\to \sigma(b)$  as $t\to b-$. Thus $g$ is HL integrable because HL integrability encloses improper integrals on finite intervals by Remarks \ref{R703.00}.
Hence (b) implies (a). 
 
 Conversely, assume that the mapping $g:[a,b]\to E$ satisfies (\ref{E30}), and is HL integrable on $[a,b]$.
We show by the Principle of Transfinite Induction, that the family $((S(\alpha)-\alpha)z_\alpha)_{\alpha\in\Lambda^{<\gamma}}$ is summable for every $\gamma\in\Lambda$. Assume that $\gamma\in \Lambda$, and that $((S(\alpha)-\alpha)z_\alpha)_{\alpha\in\Lambda^{<\beta}}$ is summable for every $\beta\in\Lambda^{<\gamma}$. If $\gamma$ is a successor, i.e., $\gamma=S(\beta)$, then  $\beta\in\Lambda^{<\gamma}$,
whence the sum $\sigma(\beta)$ of  $((S(\alpha)-\alpha)z_\alpha)_{\alpha\in\Lambda^{<\beta}}$ exists in $E$.
This result and the defining condition (s) of summability imply that $((S(\alpha)-\alpha)z_\alpha)_{\alpha\in\Lambda^{<\gamma}}$ is summable, and $\sigma(\gamma)=\sigma(\beta)+(S(\beta)-\beta)z_\beta$. 
Assume next that $\gamma$ is not a successor. Because $((S(\alpha)-\alpha)z_\alpha)_{\alpha\in\Lambda^{<\beta}}$ is summable for every $\beta\in\Lambda^{<\gamma}$, it follows from first part of this proof that for $\beta\in \Lambda^{<\gamma}$, $g$ is HL integrable integrable on $[a,\beta]$, and that (\ref{E000.4}) defines continuous mapping $f$ on  $[a,\beta]$. Thus
$$  
\sigma(\beta)=f(\beta)=\int_a^\beta g(s)\,ds, \quad \hbox{ for every} \  \beta\in \Lambda^{<\gamma}.
$$
Because $g$ is HL integrable integrable on $[a,\gamma]$, then  $\underset{\beta\to \gamma-}{\lim}\int_a^\beta g(s)\,ds$
exists and is equal to $\int_a^\gamma g(s)\,ds$ by  by Remarks \ref{R703.00}. Consequently,  $\underset{\beta\to \gamma-}{\lim}\sigma(\beta)$
exists, so that $((S(\alpha)-\alpha)z_\alpha)_{\alpha\in\Lambda^{<\gamma}}$ is summable.

The above results imply by the Principle of Transfinite Induction that the family $((S(\alpha)-\alpha)z_\alpha)_{\alpha\in\Lambda^{<\gamma}}$ is summable for every $\gamma\in\Lambda$.
In particular, $((S(\alpha)-\alpha)z_\alpha)_{\alpha\in\Lambda^{<b}}$ is summable. Thus (a) implies (b).

If (a) or (b) are valid, then both of them are  valid by the above proof. Thus the mapping $f$, defined by  (\ref{E000.4}), is a primitive of $g$, whence  $\int_a^bg(s)\,ds=f(b)-f(a)=\sigma(b)-\sigma(a)=\sigma(b)$. This proves the last conclusion.
\qed

When integrability and summability are local, we have the following result.

\begin{proposition}\label{P31} Let $\Lambda$ be a well ordered subset of a  real interval $[a,b)$, $-\infty <a < b \le\infty$,  such that
$\min\Lambda=a$ and $\sup\Lambda = b$. Assume that $g:[a,b)\to E$  is a  step mapping defined on $[a,b)$ by (\ref{E30}). 
Then the following condition are equivalent:
\newline
(a) $g$ is locally HL integrable. 
\newline
(b) The family $((S(\alpha)-\alpha)z_\alpha)_{\alpha\in\Lambda^{<b}}$ 
is locally summable.
\newline
If (a) or (b) holds, and $c\in(a,b)$, then  $\int_a^cg(t)\,dt= f(c)$, where  $f:[a,b)\to \mathbb R$ is defined on $[a,b)$ by (\ref{E000.4}).
\end{proposition}

{\bf Proof}. Assume first that the family $((S(\alpha)-\alpha)z_\alpha)_{\alpha\in\Lambda^{<b}}$ 
is locally summable. Because $\Lambda= \Lambda^{<b}$, then (\ref{E000.4}) defines
a mapping $f:[a,b)\to \mathbb R$, and $f'(t)=g(t)$ for each $t\in[a,b)\setminus\Lambda$.
As in the proof of Proposition \ref{P301} it can be shown that $f$ is continuous.                                        Thus, by Lemma \ref{L000.0}, $g$ is locally HL integrable, so that (b) implies (a). 

Conversely, assume that the mapping $g:[a,b)\to E$, defined by (\ref{E30}), is locally HL integrable on $[a,b)$.
As in the proof of Proposition \ref{P301} it can be shown that  that the family  $((S(\alpha)-\alpha)z_\alpha)_{\alpha\in\Lambda^{<\gamma}}$ is summable for every $\gamma\in\Lambda$, so that $((S(\alpha)-\alpha)z_\alpha)_{\alpha\in\Lambda}$ is locally summable. Thus (a) implies (b).

If (a) or (b) holds, then they both are valid. Assume that $c\in(a,b)$. Because   the mapping $f$, defined by (\ref{E000.4}), is a CD primitive of $g$, it follows from the last conclusion of Lemma \ref{L000.0} that $\int_a^cg(s)\,ds=f(c)-f(a)=f(c)$.
\qed   

As an application of Lemma \ref{L000.1} and Propositions \ref{P301} and \ref{P31} we obtain the following result.

\begin{proposition}\label{P32} Assume that $g:[a,\infty]\to E$ is a  step mapping satisfying  (\ref{E30}) with $b=\infty$. Then $g$ is HK integrable if and only if the family $((S(\alpha)-\alpha)z_\alpha)_{\alpha\in\Lambda^{<\infty}}$ is summable.
\end{proposition}

{\bf Proof}. Assume first that  the family $((S(\alpha)-\alpha)z_\alpha)_{\alpha\in\Lambda^{<\infty}}$ is summable.
Then it is also locally summable, whence $g$ is locally HL integrable by Proposition \ref{P32}. Thus $g$ is also locally HK integrable. 
Denote by 
 $\sigma(\gamma)$ the sum of  $((S(\alpha)-\alpha)z_\alpha)_{\alpha\in\Lambda^{<\gamma}}$, $\gamma\in\Lambda$. 
Let  $f:[a,\infty)\to E$ be defined by 
\begin{equation}\label{E040.4}
f(t)=\sigma(\gamma)+ (t-\gamma)z_\gamma, \ t\in[\gamma,S(\gamma)), \quad \gamma\in\Lambda^{<\infty}.
\end{equation}
Because the family $((S(\alpha)-\alpha)z_\alpha)_{\alpha\in\Lambda^{<\infty}}$ is summable, then
$\infty$ is a limit member of $\Lambda$. Given $\epsilon > 0$, there exists by condition (s)(ii) of summability such a $\beta_\epsilon \in\Lambda^{<\infty}$ that 
$$
\|\sigma(\beta) -\sigma(\infty)\|
<\epsilon \ \hbox{whenever $\beta\in\Lambda$ and} \ \beta_\epsilon\le \beta<\infty.
$$
If $t\in(\beta_\epsilon,\infty)$, there exists $\beta\in \Lambda$, $\beta_\epsilon\le \beta<\infty$, such that $t\in[\beta,S(\beta))$. Thus
$$
\|f(t)-\sigma(\infty)\|=\|\sigma(\beta)+(t-\beta)z_\beta-\sigma(\infty)\|<\epsilon +\|(t-\beta)z_\beta\|.
$$
Since also $\beta_\epsilon\le S(\beta)<\infty$, and since 
$$\|(t-\beta)z_\beta\|\le \|(S(\beta)-\beta)z_\beta\|=\|\sigma(S(\beta))-\sigma(\beta)\|,
$$
then
$$
\|f(t)-\sigma(\infty)\|<\epsilon+\|\sigma(S(\beta))-\sigma(\beta)\|\le\epsilon+\|\sigma({S(\beta))}-\sigma(\infty)\|+\|\sigma(\beta)-\sigma(\infty)\|<3\epsilon.
$$
This holds for every $t\in[\beta_\epsilon,\infty)$. Thus $\underset{t\to\infty}{\lim}f(t)=\sigma(\infty)$.
Because $\int_a^tg(s)\,ds=f(t)$ for each $t\in[a,\infty)$, by Proposition \ref{P32}, then
the limit $\underset{t\to\infty}{\lim}\int_a^tg(s)\,ds$ exists. This implies by Lemma \ref{L000.1} that
$g$ is HK integrable on $[a,\infty]$.

Assume next that the family $((S(\alpha)-\alpha)z_\alpha)_{\alpha\in\Lambda^{<\infty}}$ is not summable.
Then there is by  Lemma \ref{L0.00}  the greatest element $c$ in $\Lambda^{<\infty}\setminus\{\min\Lambda\}$
 such that the family $((S(\alpha)-\alpha)z_\alpha)_{\alpha\in\Lambda^{<\gamma}}$ is summable for each $\gamma\in\Lambda^{<c}$.
Moreover, $c$ is not a successor. In particular, the limit  $\underset{\gamma\to c}{\lim}\sigma(\gamma)$
does not exist. Thus $g$ is locally HL and HK integrable on $[a,c)$, but 
$\underset{\gamma\to c}{\lim} \int_a^\gamma g(s)\,ds$ does not exist. Consequently, the limit $\underset{t\to c}{\lim}\int_a^t g(s)\,ds$ does not exist, whence  Lemma \ref{L000.1} implies that $g$ is not HK integrable on $[a,c]$.  Therefore $g$ is not HK integrable on $[a,\infty]$.
\qed 

Proposition \ref{P301} is applied in the proof of the following results.

\begin{proposition}\label{P302} (a) Let $g:[a,b]\to E$,  $-\infty <a < b <\infty$,  be a step mapping that satisfies  (\ref{E30}). Then $g$
is Bochner integrable if and only if the family $((S(\alpha)-\alpha)z_\alpha)_{\alpha\in\Lambda^{<b}}$ 
is absolutely summable.
\newline
(b) Let $\Lambda$ be a well ordered subset of a  real interval $[a,b)$, $-\infty <a < b \le\infty$,  such that
$\min\Lambda=a$ and $\sup\Lambda = b$. Assume that $g:[a,b)\to E$  is a  step mapping defined on $[a,b)$ by (\ref{E30}).
Then $g$
is locally Bochner integrable if and only if the family $((S(\alpha)-\alpha)z_\alpha)_{\alpha\in\Lambda^{<b}}$ 
is locally absolutely summable.
\end{proposition}

{\bf Proof}. (a) Because $g$ is by  (\ref{E30}) strongly measurable, then $g$ is Bochner integrable if and only if the function $h=t\mapsto \|g(t)\|$ is Lebesgue integrable. Replacing $z_\alpha$ by $\|z_\alpha\|$ in (\ref{E30}) and in (\ref{E000.4}), it follows from Proposition \ref{P301} that 
$h$ is
is HL integrable if and only if the family $((S(\alpha)-\alpha)\|z_\alpha\|)_{\alpha\in\Lambda^{<b}}$ 
is summable. Because a real-valued function is HL integrable if and only if it is HK integrable, and nonnegative-valued function is HK integrable if and only if it is  Lebesgue integrable, then $h$ is Lebesgue integrable, or equivalently $g$ is Bochner integrable, if and only if the  family $((S(\alpha)-\alpha)\|z_\alpha\|)_{\alpha\in\Lambda^{<b}}$ 
is summable, or equivalently, the family $((S(\alpha)-\alpha)z_\alpha)_{\alpha\in\Lambda^{<b}}$ 
is absolutely summable.

The conclusion (b) follows from (a) and from the definitions of local integrability and local summability.  
\qed
 
\begin{proposition}\label{P303} Assume that  $\Lambda$ is a well ordered subset of a  real interval $[a,b)$ such that
$\min\Lambda=a$ and $\sup\Lambda = b$. Given a family $(z_\alpha)_{\alpha\in\Lambda}$ of $E$, let $g:[a,b]\to E$ satisfy (\ref{E30}).
\newline
(a) If $b < \infty$, then $g$ is Riemann integrable on $[a,b]$ if and only if the family $(z_\alpha)_{\alpha\in\Lambda}$
is bounded.
\newline
(b) If $b=\infty$, the improper  Riemann integral of $g$ over $[a,b]$ exists  if and only if the family $(z_\alpha)_{\alpha\in\Lambda^{<\gamma}}$ is bounded for every $\gamma\in\Lambda$, and the family $((S(\alpha)-\alpha)z_\alpha)_{\alpha\in\Lambda^{<b}}$ 
is summable.
\end{proposition}

{\bf Proof}. (a) Assume first that the family $(z_\alpha)_{\alpha\in\Lambda}$
is bounded. It follows from (\ref{E30}) that $g$ is bounded, and that, the set of its discontinuity points is a subset of $\Lambda$, and hence a null set. Thus $g$ is Riemann integrable by Lemma \ref{L000}. Conversely, if $g$ is Riemann integrable, it is bounded by Lemma \ref{L000}. Since $g(\alpha)=z_\alpha$ for each $\alpha\in \Lambda^{<b}$, then the family $(z_\alpha){\alpha\in\Lambda}$ is bounded. This proves the assertion.

(b) Assume that $b=\infty$, and that the set $(z_\alpha)_{\alpha\in\Lambda^{<\gamma}}$ is bounded for every $\gamma\in\Lambda$.
It then follows from (\ref{E30}) that $g$ is bounded, and hence Riemann integrable on each interval $[a,c]$, $a < c < \infty$. Proposition \ref{P32} implies that $g$ is HK integrable on $[a,\infty]$ if the family $((S(\alpha)-\alpha)z_\alpha)_{\alpha\in\Lambda^{<b}}$ is summable. In this case $\underset{c\to\infty}{\lim}\int_a^c g(s)\,ds$ exists by Lemma \ref{L000.1}. This limit is the improper Riemann integral because every integral $\int_a^c g(s)\,ds$, $a < c < \infty$, is Riemann integral. If the family $((S(\alpha)-\alpha)z_\alpha)_{\alpha\in\Lambda^{<b}}$ is not summable, then $g$ is not HK integrable on $[a,\infty]$, whence the improper Riemann integal over $[a,b]$ does not exist.   
\qed  

\begin{example}\label{Ex0} Let $(y_n)_{n=0}^\infty$ be a sequence in a Banach space $E$, and let $g:[0,\infty]\to E$ be defined by 
\begin{equation}\label{E300.0}
g(t)=y_n, \ t\in[n,n+1),\ n\in\mathbb N_0, \quad g(\infty)=0.
\end{equation}    
Show that
\newline
(a) $g$ is HK integrable if and only if the series $\overset{\infty}{\underset{n=0}{\sum}}y_n$ is summable;
\newline
(b) $g$ is Bochner integrable if and only if the series $\overset{\infty}{\underset{n=1}{\sum}}y_n$ is absolutely summable;
\newline
(c) The improper Riemann integral of $g$ exist if and only if the series $\overset{\infty}{\underset{n=0}{\sum}}y_n$ is summable. 
\end{example}  

{\bf Solution}.  Denoting $\alpha_n:=n$, $n\in\mathbb N_0$,
 $\Lambda=\mathbb N_0$ and $z_{\alpha_n}=y_n$, $n\in\mathbb N_0$, 
then $g$ can be rewritten as
$$
g(t)=z_{\alpha_n}, \  t\in[\alpha_n,S(\alpha_n)), \ n\in\mathbb N_0, \quad g(\infty)=0.
$$                   
The series $\overset{\infty}{\underset{n=0}{\sum}}y_n$ is summable in ordinary or absolute sense if and only if the series  $\overset{\infty}{\underset{n=0}{\sum}}(S(\alpha_n)-\alpha_n)z_{\alpha_n}$ has the same property. Moreover, 
if $\overset{\infty}{\underset{n=0}{\sum}}y_n$ is summable, then the set $\{z_{\alpha_n}|n\in\mathbb N_0\}$ is bounded. 
Thus (a) follows from Proposition \ref{P32}, (b) from Proposition \ref{P302}, and (c) from Proposition \ref{P303}.
\qed

\begin{example}\label{Ex1} Let $(y_n)_{n=0}^\infty$ be a sequence in a Banach space $E$, and let $g:[0,1]\to E$ be defined by 
\begin{equation}\label{E000.5}
g(t)=y_n,\ t\in[1-2^{-n+1},1-2^{-n}),\ n\in\mathbb N_0, \quad g(1)=0
\end{equation}    
Show that
\newline
(a) properties: $g$ is HL integrable, and the series $\overset{\infty}{\underset{n=1}{\sum}}2^{-n}y_n$ is summable, are equivalent;
\newline
(b) $g$ is Bochner integrable if and only if the series $\overset{\infty}{\underset{n=1}{\sum}}2^{-n}y_n$ is absolutely summable;
\newline
(c) The improper Riemann integral of $g$ exists if and only if the sequence $(y_n)$ is bounded.
\end{example}  

{\bf Solution}.  The correspondence $n\leftrightarrow \alpha_n:=1-2^{-n+1}$ is an order preserving isomorphism between $\mathbb N$
and $\Lambda=\{\alpha_n|n\in\mathbb N_0\}$. Denoting $z_{\alpha_n}=y_n$, $n\in\mathbb N$, and noticing that
$S(\alpha_n)-\alpha_n= 1-2^{-n}-(1-2^{-n+1})=2^{-n}$, then $g$ can be rewritten as
$$
g(t)=z_{\alpha_n},\  t\in[\alpha_n,S(\alpha_n)),\ n\in\mathbb N_0, \quad g(1)=0.
$$                   
The series $\overset{\infty}{\underset{n=1}{\sum}}2^{-n}y_n$ is summable in ordinary or absolute sense if and only if the series  $\overset{\infty}{\underset{n=1}{\sum}}(S(\alpha_n)-\alpha_n)z_{\alpha_n}$ has the same property.
Thus the conclusions of (a), (b) and (c) follow from Propositions \ref{P301}, \ref{P302} and \ref{P303}, respectively.
\qed

In view of Example \ref{Ex201} the preceding example can be generalized as follows.

\begin{example}\label{Ex302} Given $m\in \mathbb N$, let
 $\Lambda_m = \{\alpha(n_0,\dots,n_m)| \
  n_0,\dots,n_m\in\mathbb N_0\}$ be defined by (\ref{E202}). Then
$
S(\alpha(n_0,\dots,n_m))=\alpha(n_0,\dots,n_m+1)$, $ m, n_0,\dots,n_m\in\mathbb N_0$ so that
$$
S(\alpha(n_0,\dots,n_m))-\alpha(n_0,\dots,n_m)=2^{-\sum_{k=0}^mn_k-m-1}.
$$
Thus, if $(z_{\alpha(n_0,\dots,n_m)})_{\alpha(n_0,\dots,n_m)\in\Lambda_m}$ is such a family of real numbers that
 the  family 
$$
(2^{-\sum_{k=0}^mn_k-m-1}z_{\alpha(n_0,\dots,n_m)})_{\alpha(n_0,\dots,n_m)\in\Lambda_m}
$$
is summable, then the mapping $g:[0,1]\to E$, defined by
$$
g(t)=z_{\alpha(n_0,\dots,n_m)}, \ t\in[\alpha(n_0,\dots,n_m),S(\alpha(n_0,\dots,n_m))),\ \alpha(n_0,\dots,n_m)\in\Lambda_m, \quad g(1)=0,
$$
is HL integrable by Proposition \ref{P301}. According to Proposition \ref{P302} $g$ is Bochner integrable if and only if the above  family is absolutely summable. If 
the family $(z_{\alpha(n_0,\dots,n_m)})
_{\alpha(n_0,\dots,n_m)\in\Lambda_m}$, is bounded, then $g$ is Riemann integrable by Proposition \ref{P303}.
\end{example}

\begin{remark}\label{R301} 
Example \ref{Ex0} contains the results of \cite[Theorem 4 (a) and (c)]{S04}. As for related results, see \cite{VM08}. 

Let $g:[0,1)\to E$ be as in Example \ref{Ex1}, and let  $h:[0,1]\to E$ be defined by 
\begin{equation}\label{E000.6}
h(t)=y_n, \ t\in(2^{-n},2^{-n+1}], \ n\in\mathbb N_0. \quad h(0)=0.
\end{equation}    
Because $h(t)=g(1-t)$, $t\in(0,1)$, it follows from Example \ref{Ex1} that 
\newline
(a) properties: $h$ is HL integrable, and the series $\overset{\infty}{\underset{n=1}{\sum}}2^{-n}y_n$ is summable, are equivalent, 
and that
\newline
(b) $h$ is Bochner integrable if and only if the series $\overset{\infty}{\underset{n=1}{\sum}}2^{-n}y_n$ is absolutely summable.

The result (a) contains the result (a) of \cite[Proposition 5.4.1]{Sye05} and improves the results of  \cite[Proposition 5.4.2]{Sye05} and \cite[Example]{DPM02}, where  unconditional convergence of series $\overset{\infty}{\underset{n=1}{\sum}}2^{-n}y_n$ is shown to imply the HL integrability of $h$. The result (b) is equivalent to the result (c) of \cite[Proposition 5.4.1]{Sye05}.

Example \ref{Ex302} can be used to generalize further the  results of \cite{DPM02,Sye05} cited above.

In \cite{FaMa97} a notion of convergence for multiple series is defined and shown to be equivalent to the HK integrability of the associated step function  over an unbounded multidimensional interval.
\end{remark}

\section{On HK, HL, Bochner and Riemann integrability of right regulated mappings}\label{S4}
\setcounter{equation}{0}
Applying  Lemmas \ref{L400} and \ref{L401} and the results derived for step mappings in Section \ref{S3}  
we shall derive in this section criteria for HK, HL, Bochner and Riemann integrability of  right regulated mappings.

\begin{proposition}\label{P401} Given a right regulated mapping  $g:[a,b]\to E$, $-\infty < a< b<\infty$, and a positive number $\epsilon$, let $\Lambda_\epsilon$ be the well ordered subset of $[a,b]$ defined by (\ref{E403}).
Then the following properties are equivalent.
\newline 
(a) $g$ is HL integrable.
\newline
(b) The step mapping $g_\epsilon:[a,b]\to E$, defined by 
\begin{equation}\label{E405}
g_\epsilon(t)=g(\beta+),\ t\in[\beta,S(\beta)),\ \beta\in \Lambda_\epsilon^{<b}, \quad g_\epsilon(b)=g(b),
\end{equation} 
 is HL integrable.
\newline 
(c) The  family $((S(\beta)-\beta)g(\beta+))_{\Lambda_\epsilon^{<b}}$ is summable.
\end{proposition}

{\bf Proof}. It follows from Lemma \ref{L400} and (\ref{E405}) that  $\|g_\epsilon(t)-g(t)\|\le \epsilon$ whenever $t\in(\beta,S(\beta))$ and $\beta\in\Lambda_\epsilon^{<b}$. 
Because $g$ is strongly measurable by Lemma \ref{L401} and $g_\epsilon$ is strongly measurable by definition (\ref{E405}),  then $g_\epsilon-g$ is Bochner integrable, and hence also  HL integrable.  
Consequently, if $g$ is  HL integrable, then $g_\epsilon= g+(g_\epsilon-g)$ is HL integrable, and if $g_\epsilon$ is HL integrable, then $g=g_\epsilon -(g-g_\epsilon)$ is  HL integrable. This proves that
(a) and (b) are equivalent.
The equivalence of (b) and (c) follows from Proposition \ref{P301}.
\qed  

Proposition \ref{P302} and the proof of Proposition \ref{P401} is used to prove the following results.

\begin{proposition}\label{P402} Let $g:[a,b]\to E$, $-\infty < a < b < \infty$, be right regulated, and let $\epsilon$ be a positive number. 
Let $\Lambda_\epsilon$ be the well ordered subset of $[a,b]$ defined by (\ref{E403}).
Then the following properties are equivalent.
\newline 
(a) $g$ is Bochner integrable.
\newline
(b) The mapping $g_\epsilon:[a,b]\to E$, defined by (\ref{E405}),
is Bochner integrable.
\newline 
(c) The  family $((S(\beta)-\beta)g(\beta+))_{\Lambda_\epsilon^{<b}}$ is absolutely summable.
\end{proposition}

{\bf Proof}. It follows from the proof of Proposition \ref{P401} that $g_\epsilon-g$ is Bochner integrable.
Hence, if $g$ is Bochner integrable, then $g_\epsilon= g+(g_\epsilon-g)$ is Bochner integrable, and if $g_\epsilon$ is Bochner integrable, then $g=g_\epsilon -(g-g_\epsilon)$ is Bochner integrable. This implies that (a) and (b) are equivalent.
The equivalence of (b) and (c) follows from Proposition \ref{P302}.
\qed  

\begin{proposition}\label{P403} Let $g:[a,b]\to E$, $-\infty <a <b<\infty$, be right regulated,  and  
let $\Lambda_\epsilon$, $\epsilon > 0$, be the well ordered subset of $[a,b]$ defined by (\ref{E403}).
Then the following properties are equivalent.
\newline 
(a) $g$ is Riemann integrable. 
\newline
(b) $g$ is bounded.
\newline
(c) The mapping $h_\epsilon:[a,b]\to E$, defined by 
\begin{equation}\label{E4051}
h_\epsilon(t)=g(\beta+),\ t\in[\beta,S(\beta)),\ \beta\in \Lambda_\epsilon^{<b}, \quad h_\epsilon(\beta)=g(\beta), \ \beta\in\Lambda_\epsilon,
\end{equation} 
 is HL integrable.is Riemann integrable.
\newline 
(d) The families $(g(\beta))_{\beta\in\Lambda_\epsilon^{<b}}$ and $(g(\beta+))_{\beta\in\Lambda_\epsilon^{<b}}$ are bounded.
\end{proposition}

{\bf Proof}. The set of discontinuity points of $g$ is countable, whence the equivalence of (a) and (b) follows \cite[Theorem 24.45]{Sch97}. The mapping $h_\epsilon-g$ is bounded and has only a countable number of discontinuities, so that it is Riemann integrable. Hence, if $g$ is Riemann integrable, then $h_\epsilon= g+(h_\epsilon-g)$ is Riemann integrable, and if $h_\epsilon$ is Riemann integrable, then $g=h_\epsilon -(g-h_\epsilon)$ is Riemann integrable. This implies that (a) and (c) are equivalent.
The proof of the equivalence of (c) and (d) is similar to that of Proposition \ref{P303}.
\qed

Now we are in position to prove the results presented in the Introduction.
\begin{theorem}\label{T41} (The Fundamental Theorem of Calculus for right regulated mappings) Assume that a mapping
 $g:I\to E$, $-\infty<\min I<\sup I\le\infty$ is right regulated.
\begin{enumerate}
\item [(a)] $g$ is locally HL integrable if and only if it has a CD primitive.
\item [(b)] $g$ is locally Bochner integrable if and only if it has a locally absolutely continuous CD primitive.
\item [(c)] $g$ is locally Riemann integrable if and only if it is locally bounded, in which case $g$ has a locally Lipschitz continuous CD primitive.
\end{enumerate}
\end{theorem}
{\bf Proof}. (a) Denote $a=\min I$ and $c=\sup I$. It follows from Lemma \ref{L000.0} that
if $g$ has a CD primitive, then $g$ is locally HL integrable. 
To prove converse, assume that $g$ is locally HL integrable. Given $b\in (a,c)$,
define for each $n\in\mathbb N$ the step mapping $g_n:[a,b]\to E$ by 
\begin{equation}\label{E406}
g_n(t)=g(\beta+),\ t\in[\beta,S(\beta)),\ \beta\in \Lambda_{\frac 1n}^{<b}, \quad g_n(b)=g(b).
\end{equation}
Because $g$ is HL integrable on $[a,b]$, it follows from 
Proposition \ref{P401} that the family $((S(\beta)-\beta)g_n(\beta+))_{\Lambda_{\frac 1n}^{<b}}$ is summable
for every $n\in\mathbb N$. Denote by $\sigma_n(\gamma)$ the sum of  the family 
$((S(\beta)-\beta)g(\beta+))_{\Lambda_{\frac 1n}^{<\gamma}}$, $\gamma\in \Lambda_{\frac 1n}^{<b}$. The proof of
Proposition \ref{P301} that implies that for each $n\in\mathbb N$ 
the mapping $f_n:[a,b)\to E$, defined by
\begin{equation}\label{E400.4}
f_n(t)=\sigma_n(\gamma)+ (t-\gamma)g(\gamma+), \ t\in[\gamma,S(\gamma)), \ \gamma\in\Lambda_{\frac 1n}^{<b},
\end{equation}
is a CD primitive of $g_n$. Thus, for each $n\in \mathbb N$,  the mapping $f_n$ is continuous, $f_n'(t)=g_n(t)$ for all $t\in[a,b]\setminus \Lambda_{\frac 1n}$, and $f_n(a)=\sigma_n(a)=0$. Moreover, if $g(\beta+)=g(\beta)$ for $\beta\in Z=\bigcup_{n=1}^\infty \Lambda_n$, then $\|g_n(t)-g(t)\|\le \frac 1n$ for each $t\in[a,b)$ by Lemma \ref{L400} and (\ref{E406}), so that the sequence $(g_n)_{n=1}^\infty$
converges uniformly to $g$. Consequently, it follows from \cite[(8.6.4)]{Die60} that
the sequence $(f_n)_{n=1}^\infty$ converges uniformly on $[a,b)$ to a continuous mapping $f:[a,b)\to E$, and
$f'(t)=g(t)$ for each $t\in [a,b)\setminus Z$. 
$f$ has these properties also when right continuity of $g$ in $Z$ is not assumed. Because $Z$ is countable, then $f$ is a CD primitive of the restriction of $g$ to $[a,b]$.      

Choose an increasing sequence $(c_n)_{n=1}^\infty$ from $(a,c)$  
so that it converges to $c$. The interval $[a,c)$  is the union of increasing sequence of compact intervals $[a,c_n]$, and $g$ is HL integrable on these compact intervals. By the above proof the restriction of $g$ to $[a,c_n]$ has a CD primitive $f_n$, and $f_n(a)=0$ for each $n\in\mathbb N$. Defining
$$
f(t)=\begin{cases} f_1(t)-f(a), \ &t\in[a,c_1),\\ f_{n+1}(t)-f_{n+1}(a), \ &t\in[c_n,c_{n+1}), \ n\in\mathbb N,\end{cases}
$$ 
we obtain a mapping $f:[a,c)\to E$ which is a CD primitive of $g$ (cf. Remark after \cite[(8.7.1)]{Die60}). 

(b) If $g$ has a locally absolutely continuous CD primitive $f$, then $g$ is locally Bochner integrable by \cite[Theorem 1.4.6]{HL94}. Conversely, assume that $g$ is locally Bochner integrable, and let $[a,b]$ be a compact subinterval of $I$. Then $g$ is Bochner integrable on $[a,b]$, whence the restriction of $g$ to $[a,b]$ has  by \cite[Theorem 1.4.6]{HL94} an absolutely continuous primitive $h:[a,b]\to E$.
 $g$ is also locally HL integrable by  \cite[Proposition 3.6.3 and Theorem 5.1.4]{Sye05}. Thus $g$ has by the proof of (a) a CD primitive $f:I\setminus\{\sup I\}\to E$.
It follows from Lemma \ref{L000.0} and from the definition (\ref{E200}) of the Henstock-Kurzweil integral that
$
 f(x)-f(a)=h(x)-h(a)$, i.e., 
$f(x)= h(x)+f(a)-h(a)$ for all $x\in[a,b]$. Thus $f$ is absolutely continuous on $[a,b]$.
Consequently, $f$ is an absolutely continuous CD primitive of the restriction of $g$ to $[a,b]$.     

(c) Assume that $g$ is locally Riemann integrable, and  let $[a,b]$ be a compact subinterval of $I$. It follows from Proposition \ref{P403} that $g$ is Riemann integrable on $[a,b]$ if and only if $g$ is bounded on $[a,b]$, in which case there is such a positive constant $M$ that
$\|g(t)\|\le M$ for each $t\in[a,b]$. Because $g$ is also locally HL integrable, it has a  CD primitive $f$ by the proof of (a),
and $f(b)-f(a)=\int_a^bg(t)\,dt$ by Lemma \ref{L000.0}. Thus,  $\|f(b)-f(a)\|\le\int_a^b\|g(t)\,dt\|\le M(b-a)$.
This holds for every compact subinterval $[a,b]$ of $I$, whence $f$ is locally Lipschitz continuous.
\qed

The following results are easy consequences of the results of Theorem \ref{T41} and Lemma \ref{L000.1} and the definitions of integrals and improper integrals. 

\begin{corollary}\label{C41} Let $g:I\to E$, $-\infty<\min I<\sup I\le\infty$ be right regulated.
\begin{enumerate}
\item [(a)] $g$ is HL integrable when $I$ is bounded if and only if $g$ has a CD primitive that has the left limit at $\sup I$.
\item [(b)] $g$ is HK integrable if  it has a  CD primitive that has the left limit at $\sup I$.
\item [(c)] $g$ is Bochner integrable if and only if the function it has a locally absolutely continuous CD primitive that has the left limit at $\sup I$.
\item [(d)] $g$ is Riemann integrable if and only if it is bounded and $I$ is bounded.
\item [(e)] The improper Riemann integral of $g$ from $\min I$ to $\sup I$  if $g$ is locally bounded, and its CD primitive has the left limit at $\sup I$.  
\end{enumerate}
\end{corollary}
The next result follows from Lemma \ref{L0.00} and Propositions \ref{P401}, \ref{P402} and \ref{P403}. 

\begin{theorem}\label{T42} Let $g:I\to E$, $-\infty<\min I<\sup I\le\infty$ be right regulated. 
\begin{enumerate}
\item  [(a)] For each compact subinterval $[a,b]$ of $I$, either $g$ is Riemann integrable on $[a,b]$, or there  exists the greatest number $c_1$ in $(a,b]$ such that $g$ is locally Riemann integrable on $[a,c_1)$
\item  [(b)] For each compact subinterval $[a,b]$ of $I$, either $g$ is Bochner integrable on $[a,b]$, or there  exists the greatest number $c_2$ in $(a,b]$ such that $g$ is locally Bochner integrable on $[a,c_2)$.
\item  [(c)] For each compact subinterval $[a,b]$ of $I$, either $g$ is HL integrable on $[a,b]$, or there  exists the greatest number $c_3$ in $(a,b]$ such that $g$ is 
locally HL integrable on $[a,c_3)$. 
\end{enumerate}
\end{theorem}

{\bf Proof}. Let $[a,b]$ be a compact subinterval of $I$, let $\epsilon$ be a positive number, and 
let $\Lambda_\epsilon$ be the well ordered subset of $[a,b]$ defined by (\ref{E403}).
\newline
(a) According to Lemma \ref{L0.00} (a) the family $(g(\beta+))_{\beta\in\Lambda_\epsilon^{<b}}$ is bounded, or there exists
the greatest number $c_1$ in $\Lambda_\epsilon^{<b}$, $c_1>a$, such that  the family $(g(\beta+))_{\beta\in\Lambda_\epsilon^{<\gamma}}$ is bounded for every $\gamma\in\Lambda_\epsilon^{<c_1}$. 
This result and Proposition \ref{P403} imply that $g$ is Riemann integrable either on $[a,b]$, or on 
$[a,\gamma]$, for every $\gamma\in\Lambda_\epsilon^{<c_1}$. This proves (a) because $c_1$ is by Lemma \ref{L0.00} (d)  not a successor.   

(b) By Lemma \ref{L0.00} (b) the family $((S(\beta)-\beta)g(\beta+))_{\beta\in\Lambda_\epsilon^{<b}}$ is absolutely summable, or there exists
the greatest number $c_2$ in $\Lambda_\epsilon^{<b}$, $c_1>a$, such that  the family $((S(\beta)-\beta)g(\beta+))_{\beta\in\Lambda_\epsilon^{<\gamma}}$ is absolutely summable for every $\gamma\in\Lambda_\epsilon^{<c_2}$. 
This result implies by Proposition \ref{P402} that $g$ is Bochner integrable either on $[a,b]$, or on 
$[a,\gamma]$, for every $\gamma\in\Lambda_\epsilon^{<c_2}$. This implies conclusion (b), since by Lemma \ref{L0.00} (d) $c_2$ is not a successor. 

(c) The proof of (c) is similar to that of (b) when  absolute summability is replaced by summability, Lemma \ref{L0.00} (b) by Lemma \ref{L0.00} (c), and   Proposition \ref{P402} by  Proposition \ref{P401}.
\qed 

\begin{example}\label{Ex401} Denote $c_0=\{x=(x_i)_{i=1}^\infty|x_i\in\mathbb R,\,i\in\mathbb N, \ \underset{i\to\infty}{\lim}x_i=0\}$. $c_0$ is a vector space with respect to componentwise addidion and scalar multiplication, and $\|x\|=\underset{i\in\mathbb N}{\sup}|x_i|$ defines a Banach norm in $c_0$.  Define a mapping $g_0:\mathbb R_+\to c_0$ by 
\begin{equation}\label{E47}
g_0(t)=\left(\sum_{n=1}^\infty\frac 1{n^2i}\left(2(nt-\left\lfloor nt\right\rfloor)\cos\left(\frac {\pi}{2(nt-\left\lfloor nt\right\rfloor)}\right)+\frac {\pi}{2}\sin\left(\frac{\pi}{2(nt-\left\lfloor nt\right\rfloor)}\right)\right)\right)_{i=1}^\infty, \ t\in\mathbb R_+, 
\end{equation}
where $\left\lfloor nt \right\rfloor=m,\ m-1<nt\le m, \ m=0,1,\dots$.  $g_0$ is right regulated.
The set $\mathbb Q_+$ of all rational numbers of $\mathbb R_+$ is the set of discontinuity points $g_0$ (cf. \cite[(236)]{EL40}). Moreover, all these discontinuities are of second kind. A CD primitive of $g$ is given by
\begin{equation}\label{E48}
f(t)=\left(\sum_{n=1}^\infty\frac {(nt-\left\lfloor nt\right\rfloor)^2}{n^3i}\cos\left(\frac \pi{2(nt-\left\lfloor nt\right\rfloor)}\right)\right)_{i=1}^\infty, \ t\in\mathbb R_+. 
\end{equation}
Because $g$ is bounded, it is also locally Riemann integrable by Theorem \ref{T41}. 

The mapping $g=t\mapsto e^{-t}g_0(t)$ has the improper Riemann integral $\int_0^\infty g(t)\,dt$.
\end{example}

\begin{example}\label{Ex402} Let $g_0$ and $f_0$ be defined by (\ref{E47}) and (\ref{E48}). Define mappings $g_m:\mathbb R_+\to c_0$, $m\in\mathbb N$, by  
\begin{equation}\label{E407}
g_m(t)=g_0(t)+\left(\frac{1}{i}\sum_{n=1}^{i\wedge m}\left(\cos\left(\frac{\pi}{2(nt-\left\lfloor nt\right\rfloor)}\right) +\frac{\pi\sin(\frac{\pi}{2(nt-\left\lfloor nt\right\rfloor)})}{2(nt-\left\lfloor nt\right\rfloor)}\right)\right)_{i=1}^\infty, \ t\in\mathbb R_+, 
\end{equation}
where $i\wedge m=\min\{i,m\}$. $g_m$ is right regulated, and $\mathbb Q_+$ is the set of discontinuity points, of second kind, of $g_m$. for all $m\in\mathbb N$. The mapping $f_m:\mathbb R_+\to c_0$, defined by
\begin{equation}\label{E408}
f_m(t)=f_0(t)+\left(\sum_{n=i}^{i\wedge  m}\frac{(nt-\left\lfloor nt\right\rfloor)^2}{ni}\cos\left(\frac{\pi}{2(nt-\left\lfloor nt\right\rfloor)}\right)\right)_{i=1}^\infty, \ t\in\mathbb R_+, 
\end{equation}
is a CD primitive of $g_m$ for each $m\in\mathbb N$.  It then follows from Theorem \ref{T41} that the mappings $g_m$ are locally HL integrable. On the other hand, $g_m$ is neither locally  Bochner integrable nor locally Riemann integrable for any $m\in\mathbb N$.

The mappings $t\mapsto e^{-t}g_m(t)$ are HK integrable on $\mathbb R_+$.
\end{example}

\begin{example}\label{Ex403} Let $g_0$ and $f_0$ be defined by (\ref{E47}) and (\ref{E48}).  Define mappings $g^m:\mathbb R_+\to c_0$, $m\in\mathbb N$,  by  
\begin{equation}\label{E409}
g^m(t)=g_0(t)+\left(\frac 1i\sum_{n=1}^{i\wedge  m}\frac{1}{2\sqrt{\left\lfloor nt\right\rfloor-nt}}\right)_{i=1}^\infty, \ t\in\mathbb R_+. 
\end{equation}
$g^m$ is right regulated, and $\mathbb Q_+$ is its set of discontinuity points, of second kind, for every $m\in\mathbb N$. The mappings $f^m:\mathbb R_+\to c_0$, defined by
\begin{equation}\label{E410}
f^m(t)=f_0(t)+\left(\sum_{n=1}^{i\wedge m}\frac{\left\lfloor nt\right\rfloor-\sqrt{\left\lfloor nt\right\rfloor-nt}}{ni}\right)_{i=1}^\infty, \ t\in\mathbb R_+, 
\end{equation}
are absolutely continuous, and $(f^m)'(t)=g^m(t)$ for all $t\in\mathbb R_+\setminus\mathbb Q_+$.
Hence, every $g^m$ is locally Bochner integrable by Theorem \ref{T41}. But $g^m$ is not locally bounded, and hence not locally Riemann integrable, for any $m\in\mathbb N$.

The mappings $t\mapsto e^{-t}g^m(t)$ are Bochner integrable on $\mathbb R_+$.
\end{example}

\begin{remark}\label{R401} 
Integrability results derived in  Propositions \ref{P401}, \ref{P402} and \ref{P403}, and in Theorems \ref{T41} and \ref{T42} for right regulated mappings have  also analogous counterparts for left regulated mappings.
\end{remark}

\section{Applications to impulsive differential equations}\label{5}
\setcounter{equation}{0}
Let $E$ be a Banach space and $[a,c)$, $-\infty < a < c\le\infty$, a real interval. Denote by $HL_{loc}([a,c),E)$ the space of all locally HL integrable mappings from $[a,c)$ to $E$. Almost everywhere (a.e.) equal mappings of $HL_{loc}([a,c),E)$ are identified.  Consider the following
impulsive problem
\begin{equation}\label{E512.46}
u'(t)=f(t,u) \ \hbox{ a.e. on } \ [a,c), \quad \Delta u(\lambda)=D(\lambda,u), \ \lambda\in \Lambda,
 \end{equation}
where  $f:[a,b)\times HL_{loc}([a,c),E)\to E$,
$\Delta u(\lambda)=u(\lambda+)-u(\lambda)$, $D:\Lambda\times
HL_{loc}([a,c),E)\to E$, and $\Lambda$ is a well ordered 
subset of $[a,c)$ with $a=\min\Lambda$ and $c=\sup\Lambda$. When $t\in [a,c)$, we denote $\Lambda^{<t}=\{\lambda\in \Lambda: \lambda < t\}$. If a  family $(x(\lambda))_{\lambda\in \Lambda}$ of $E$ is locally summable, and $t\in[a,c)$,
denote   by $\underset{{\lambda\in\Lambda^{<t}}}{\sum}x(\lambda)$ the sum of the family $(x(\lambda))_{\lambda\in \Lambda^{<t}}$.

We say that $u:[a,c)\to E$ is a solution
of problem (\ref{E512.46}) if it satisfies the equations of
(\ref{E512.46}), and if it belongs to the set
$$
V=\{u\in HL_{loc}([a,c),E)| u\ \mbox{ is a.e. differentiable and right continuous}\}.
$$
The following result  allows us to transform problem (\ref{E512.46}) into an
integral equation.

\begin{lemma}\label{L512.46}
Let $g\in HL_{loc}([a,c),E)$ and assume that 
a  family $(z(\lambda))_{\lambda\in \Lambda}$ of $E$ is locally summable. Then the problem
\begin{equation}\label{E512.47}
u'(t)=g(t) \ \hbox{ a.e. on } \ [a,c),\quad \Delta u(\lambda)=z(\lambda), \ \lambda\in \Lambda,
 \end{equation}
has a unique solution $u$.  
 This solution can be represented as
\begin{equation}\label{E512.48}
u(t)= \sum_{\lambda\in\Lambda^{<t}}z(\lambda) +
\int_a^tg(s)ds, \quad  t\in [a,c).
\end{equation}
Moreover, $u$ is increasing with respect to $g$ and $z$.
\end{lemma}

\noindent{\bf Proof:}
Let $u:[a,c)\to E$ be defined by (\ref{E512.48}).
It is easy to verify that
\begin{equation}\label{E512.49}
u'(t)=g(t) \ \hbox{ for a.e. } \ t\in [a,c).\end{equation}
For each $\lambda\in \Lambda$ the open interval $(\lambda,S(\lambda))$  does
not contain any point of $\Lambda$, so that
\begin{equation}\label{E512.410}
\Delta u(\lambda)=u(\lambda+)-u(\lambda) = \lim_{t\to
\lambda+}\left(z(\lambda)+\int_\lambda^tg(s)ds\right)=z(\lambda), \ \lambda\in \Lambda.
\end{equation}
It follows from (\ref{E512.48})  that
\begin{equation}\label{E512.411}
u(t)= u_1(t) + u_2(t),
\end{equation}
where
$$u_1(t)= \int_{a}^tg(s)ds, \quad u_2(t)=\sum_{\lambda\in\Lambda^{<t}}z(\lambda),
 \quad t\in [a,c).
$$
Because $(z(\lambda))_{\lambda\in \Lambda}$ is locally summable, then both $u_1$ and $u_2$ belong to
$V$.
This, (\ref{E512.49}) and (\ref{E512.410})  imply
that $u$ is a solution of problem (\ref{E512.47}).

If $v\in V$ is a solution of (\ref{E512.47}), then $w=u-v$ belongs to
$V$, and $\Delta w(\lambda)=0$ for each $\lambda\in \Lambda$,
whence $w$ is a solution of the initial
value of problem
\begin{equation}\label{E512.4111}
w'(t)=0 \ \hbox{ a.e. on } \ [a,c), \quad w(a)=0.
\end{equation}
This  implies that $w(t)\equiv 0$,
i.e., $u=v$.

The last assertion of Lemma is a direct consequence from the
representation (\ref{E512.48}) and \cite[Lemma 9.11]{CH11}. \qed

Assume that $g\in HL_{loc}([a,c),E)$ is right regulated. 
Given $\epsilon > 0$ and $b\in(a,c)$, let $\Lambda_\epsilon$ be the well ordered subset of $[a,b]$ defined by (\ref{E403}), and let $g_\epsilon:[a,b]\to E$ be defined by (\ref{E405}).
Because $g$ is locally HL integrable, it follows from Proposition \ref{P401} that 
$g_\epsilon$ is  HL integrable on $[a,b]$, and that  
 the family $((S(\alpha)-\alpha)g(\alpha+))_{\alpha\in\Lambda_\epsilon^{<b}}$ is summable.
Let $\sigma_\epsilon(\gamma)$ denote the sum of the family $((S(\alpha)-\alpha)g(\alpha+))_{\alpha\in\Lambda_\epsilon^{<\gamma}}$, $\gamma\in\Lambda_\epsilon$. 
Define a mapping $f_\epsilon:[a,c)\to E$ by
\begin{equation}\label{E511.4}
f_\epsilon(t)=\sigma_\epsilon(\gamma)+ (t-\gamma)g(\gamma+), \ t\in[\gamma,S(\gamma)), \ \gamma\in\Lambda_\epsilon^{<c}.
\end{equation}
By the proof of Proposition \ref{P301}, $f_\epsilon$ is a CD primitive of $g_\epsilon$.
It follows from Lemma \ref{L400} that $\|g_\epsilon(t)-g(t)\|\le\epsilon$ for all $t\in [a,b]\setminus \Lambda_\epsilon$.
Thus
$$
\left\|\int_a^tg(s)\,ds-\int_a^tg_\epsilon(s)\,ds\right\|= \left\|\int_a^tg(s)\,ds-f_\epsilon(t)\right\|\le\epsilon(t-a).
$$
The above considerations and Theorem \ref{T41} imply the following results for solutions of problem 
(\ref{E512.47}).
\begin{proposition}\label{P511.01}  Assume that $g\in HL_{loc}([a,c),E)$ is right regulated, and that the  family $(z(\lambda))_{\lambda\in \Lambda}$ is summable. Then for all fixed $b\in(a,c)$ and $\epsilon > 0$ the mapping $u_\epsilon:[a,b]\to E$, defined by 
\begin{equation}\label{E512.50}
u_\epsilon(t)= \sum_{\lambda\in\Lambda^{<t}}z(\lambda) + f_\epsilon(t), \quad  t\in [a,b],
\end{equation}
approximates the solution of problem (\ref{E512.47}) uniformly on $[a,b]$ within the accuracy $\epsilon(b-a)$.
The differential equation of (\ref{E512.47}) holds in the complement of a countable subset of $[a,b]$.
\end{proposition}

In what follows we assume that $E$ is a Banach space ordered by a regular order cone,  and that  the function space $HL([a,b],E)$ is ordered a.e. pointwise. The following fixed point result
is a consequence of \cite[Theorem 2.16 and Proposition 9.39]{CH11}.

\begin{theorem}\label{T511.00}
Let $[w_-,w_+]=\{g\in HL_{loc}([a,c),E)|w_-\le g\le w_+\}$ be a nonempty
order interval in $HL_{loc}([a,c),E)$. Then every increasing mapping
$G:HL_{loc}([a,b),E)\to[w_-,w_+]$ has the smallest and greatest
fixed points, and they are increasing with respect to $G$.
\end{theorem}

Let us impose the following hypotheses on the mappings $f$ and $D$ in problem
(\ref{E512.46}).
\begin{description}
\item[(f0)] There exist locally HL integrable mappings $f_\pm:[a,c)\to E$ such that $f_-(t)\le f(t,u)\le f_+(t)$ for a.e. $t\in [a,c)$ and for all $u\in HL_{loc}([a,c),E)$.
\item[(f1)] The mapping $f(\cdot,u)$ is right regulated for each $u\in HL_{loc}([a,c),E)$.
\item[(f2)] $f(s,\cdot)$ is increasing for a.e. $s\in [a,c)$.
\item[(D0)] $D(\lambda,\cdot)$ is increasing for all $\lambda\in \Lambda$, and
there exist $z_\pm:\Lambda\to E$ such that
$z_-(\lambda)\le D(\lambda,u)\le z_+(\lambda)$ for all $\lambda\in
\Lambda$ and $u\in HL_{loc}([a,c),E)$, and that the families
$(z_\pm(\lambda))_{\lambda\in \Lambda}$ are locally summable.
\end{description}

As an application of Theorem \ref{T511.00} we get the following
existence and comparison result for problem (\ref{E512.46}).

\begin{theorem}\label{T512.44}
Let the mappings $f$ and $D$  in (\ref{E512.46}) satisfy the
hypotheses (f0)--(f2)  and (D0). Then problem (\ref{E512.46}) has the
smallest and greatest solutions $u_*$ and $u^*$ in $V$. Moreover,
these solutions are increasing with respect to $D$ and $f$, and they satisfy the differential equation of (\ref{E512.46}) the complement of a countable subset of $[a,c)$.
\end{theorem}

\noindent{\bf Proof:} The hypotheses (f0) and (D0) ensure that the equations
\begin{equation}\label{E512.412}
w_\pm(t)= \sum_{\lambda\in
\Lambda^{<t}}z_\pm(\lambda) + \int_a^tf_\pm(s)ds,
\end{equation}%
define mappings $w_\pm\in HL_{loc}([a,c),E)$.
By using the hypotheses, and \cite[Lemma 9.11 and Proposition 9.14]{CH11} it can be shown that the equation
\begin{equation}\label{E512.413}
Gu(t):= \sum_{\lambda\in
\Lambda^{<t}}D(\lambda,u) + \int_a^tf(s,u)ds,\quad t\in[a,c),
\end{equation}
defines an increasing mapping $G:HL_{loc}([a,c),E)\to [w-,w+]$. From Theorem \ref{T511.00} it then follows that
$G$ has the smallest and greatest fixed points
$u_*$ and $u^*$, and they are increasing with respect to  $D$ and
$f$. Because by Lemma \ref{L512.46} the solutions of problem
(\ref{E512.46}) are the same as the fixed points of $G$, then $u_*$ and $u^*$ are the smallest and
greatest solutions of problem (\ref{E512.46}), and they are
increasing with respect to $D$ and $f$. 
To show the validity of the last conclusion, let $u$ be any fixed point of $G$, i.e., 
 \begin{equation}\label{E512.423}
u(t)= \sum_{\lambda\in
\Lambda^{<t}}D(\lambda,u) + \int_a^tf(s,u)ds,\quad t\in[a,c).
\end{equation}
The mapping $f(\cdot,u)$ is by the hypothesis (f1) right regulated, and also locally HL integrable on $[a,c)$. 
Thus it has by Theorem \ref{T41} a CD primitive $\tilde f$, and $\int_a^tf(s,u)ds=\tilde f(t)-\tilde f(a)$, $t\in[a,c)$. Hence there is a countable subset $Z_1$ of $[a,c)$ such that $\tilde f'(t)=f(t,u)$ for each $t\in[a,c)\setminus Z_1$. Denoting $Z=Z_1\cup\Lambda$, it then follows from (\ref{E512.423}) that
$$
u'(t)= \frac d{dt}\left(\sum_{\lambda\in
\Lambda^{<t}}D(\lambda,u) + \int_a^tf(s,u)ds\right)=f(t,u), \quad t\in [a,c)\setminus Z.
$$
This proves the last conclusion.     
\qed

\begin{example}\label{ex406.502} The cone  of those elements of
$E=c_0$ with nonnegative coordinates is regular. Choose $[a,c)=[0,\infty)=\mathbb R_+$. Let $g_0:\mathbb R_+\to c_0$ be defined by (\ref{E47}), 
and define  $q_i: \mathbb
R\to \mathbb R$, $i=1,2,\dots$, by
$$
q_i(s)=\frac 1{2^i}\sum_{m=1}^i\sum_{k=1}^\infty\frac{\frac{\pi}2
+\arctan(k^{\frac 1m}s)}{(km)^2}, \ s\in \mathbb R, \ i=1,2,\dots .
$$
 For $x=(x_1,x_2,\dots)\in c_0$,
define
$$
g(t,x) = g_0(t) +
\left(q_i\Bigl(\sum_{j=1}^i x_j\Bigr)\right)_{i=1}^\infty, \quad t\in\mathbb R_+.
$$
Then one can easily verify  that $f(t,u)=g(t,u(t))$
satisfies hypotheses (f0)--(f2).

Let $\Lambda$ be a well ordered subset of real numbers with $\min\Lambda=0$ and $\sup\Lambda=\infty$. 
Denoting
$$
c(\lambda)=(c_1(\lambda), c_2(\lambda),\dots),  \hbox{ where }
c_i(\lambda)=2^{-i}z_\lambda, \ \  \lambda\in \Lambda. \
i=1,2,\dots,
$$
Assuming that the family $\underset{\lambda\in\Lambda}{\sum}z_\lambda$ is a summable family of real numbers $z_\lambda$, then 
the family $\underset{\lambda\in\Lambda}{\sum}c(\lambda)$ is summable in $c_0$. Thus the mapping $D(\cdot,u)\equiv c$ has
the properties assumed in (D0). With $c$ and $g$ defined before,
consider the problem
\begin{equation}\label{E406.506}
u'(t)= g(t,u(t)) \ \hbox { a.e. on } \ \mathbb R_+,
\quad \Delta u(\lambda)= c(\lambda), \ \lambda\in \Lambda.
\end{equation}
The above proof shows that the hypotheses  of Theorem \ref{T512.44}
are valid, when $f(t,u)=g(t,u(t))$ and $D(\lambda,u)=c(\lambda)$.
Thus problem (\ref{E406.506}) has the smallest and greatest solutions.
\end{example}

{\bf Competing interests}

The author declares that he has no competing interests.

{\bf Author's contributions}

The work is realized by the author.

"...the infinite series can be totaled at any given point, and this total (more properly, a subtotal) provides the fullness of the sweetness of goal attainment for a given person at a given time and status. But sooner or later, this same person begins to hunger and yearn for new and greater goals, and such adventures in growth will be forever forthcoming in the fullness of time and the cycles of eternity."([The Book of Urantia, 118,0,11]).
\end{document}